\documentclass[fleqn]{article}

\usepackage[a4paper, left=2cm, right=2cm, top=3cm, bottom=3cm]{geometry}                           % page format
\usepackage[applemac]{inputenc}                                                                    % Mac OS X character set
\usepackage[T1]{fontenc}                                                                           % character set
\usepackage[colorlinks, linkcolor=black, citecolor=black, urlcolor=black]{hyperref}                % hyperlinked pdf-output
\usepackage{amsmath}                                                                               % AMS-LaTeX
\usepackage{amssymb}                                                                               % AMS-LaTeX
\usepackage{amsthm}                                                                                % AMS-LaTeX
\usepackage{upgreek}                                                                               % upright greek letters
\usepackage{tikz}                                                                                  % TikZ for graphics

\usetikzlibrary{matrix}                                               

\swapnumbers
\theoremstyle{definition}
\newtheorem*{example*}{Example}
\newtheorem*{proposition*}{Proposition}

\newcommand{\booktitle}[1]{\textsl{#1}}                                                            % booknames in slanted letters
\newcommand{\eigenname}[1]{\textsc{#1}}                                                            % names in capital letters
\newcommand{\newnotion}[1]{\textit{#1}}                                                            % new notions in italic letters

\newenvironment{smallpmatrix}{\left( \begin{smallmatrix}}{\end{smallmatrix} \right)}               % small matrices with round brackets

\DeclareMathOperator{\HomotopyCategory}{\mathrm{Ho}}                                               % homotopy category
\DeclareMathOperator{\Ob}{\mathrm{Ob}}                                                             % objects
\newcommand{\Bif}{\mathbf{Bif}}                                                                    % full subcategory of bifibrant objects
\newcommand{\categoricalequivalent}{\simeq}                                                        % categorical equivalent
\newcommand{\Cof}{\mathbf{Cof}}                                                                    % full subcategory of cofibrant objects
\newcommand{\cylinderequivalence}{\mathrm{s}}	                                                      % cylinder equivalence
\newcommand{\cylinderhomotopic}{\mathrel{\stackrel{\mathrm{c}}{\sim}}}                             % cylinder homotopic
\newcommand{\congruent}{\equiv}                                                                    % congruent
\newcommand{\directsum}{\oplus}                                                                    % direct sum
\newcommand{\emb}{\mathrm{emb}}                                                                    % embedding
\newcommand{\fgMod}{\mathbf{mod}}                                                                  % category of finitely generated modules
\newcommand{\Fib}{\mathbf{Fib}}                                                                    % full subcategory of bifibrant objects
\newcommand{\homotopic}{\sim}                                                                      % homotopic
\newcommand{\ins}{\mathrm{ins}}                                                                    % insertion
\newcommand{\Integers}{\mathbb{Z}}                                                                 % set of integers
\newcommand{\isomorphic}{\cong}                                                                    % isomorphic
\newcommand{\LocalisationFunctor}[1][]{\Upgamma^{#1}}                                              % localisation functor
\newcommand{\map}{\rightarrow}                                                                     % map
\newcommand{\Naturals}{\mathbb{N}}                                                                 % set of natural numbers
\newcommand{\pathhomotopic}{\mathrel{\stackrel{\mathrm{p}}{\sim}}}                                 % path homotopic
\newcommand{\smallcoprod}{\amalg}                                                                  % small coproduct sign
\newcommand{\terminalobject}[1][]{{*}^{#1}}                                                        % terminal object
\newcommand{\weaklyequivalent}{\approx}                                                            % weakly equivalent

\tikzset{diagram/.style={matrix of math nodes, row sep=#1, column sep=#1, text height=1.6ex, text depth=0.45ex, inner sep=0pt, nodes={inner sep=0.3333em}}, diagram/.default=2.5em}
\tikzset{equality/.style={-, double}}

\title{Faithfulness of a functor of Quillen}
\author{William G.\ \\ Dwyer \and Andrei \\ R\u{a}dulescu-Banu \and Sebastian \\ Thomas}
\date{October 6, 2009}

\begin{document}

\maketitle

% \tableofcontents

\renewcommand{\thefootnote}{\fnsymbol{footnote}}
\renewcommand{\thefootnote}{\arabic{footnote}}

\begin{abstract}
There exists a canonical functor from the category of fibrant objects of a model category modulo cylinder homotopy to its homotopy category. We show that this functor is faithful under certain conditions, but not in general.
\end{abstract}

\section{Introduction} \label{sec:introduction}

We let \(\mathcal{M}\) be a model category. \eigenname{Quillen} defines in \cite[ch.\ I, \S 1]{quillen:1967:homotopical_algebra} a homotopy relation on the full subcategory \(\Fib(\mathcal{M})\) of fibrant objects, using cylinders. He obtains a quotient category \(\Fib(\mathcal{M}) / {\cylinderhomotopic}\) and a canonical functor
\[\Fib(\mathcal{M}) / {\cylinderhomotopic} \map \HomotopyCategory \Fib(\mathcal{M}).\]
The question occurs whether this functor is faithful.

We show that it is faithful if \(\mathcal{M}\) is left proper and fulfills an additional technical condition. Moreover, we show by an example that it is not faithful in general.

\subsection*{Conventions and notations}

\begin{itemize}
\item The composite of morphisms \(f\colon X \map Y\) and \(g\colon Y \map Z\) is denoted by \(f g\colon X \map Z\).
\item Given \(n \in \Naturals_0\), we abbreviate \(\Integers / n := \Integers / n \Integers\). Given \(k, m, n \in \Naturals_0\), we write \(k\colon \Integers / m \map \Integers / n, a + m \Integers \mapsto k a + n \Integers\), provided \(n\) divides \(k m\).
\item Given a category \(\mathcal{C}\) with finite coproducts and objects \(X, Y \in \Ob \mathcal{C}\), we denote by \(X \smallcoprod Y\) a (chosen) coproduct. The embedding \(X \map X \smallcoprod Y\) is denoted by \(\emb_0\), the embedding \(Y \map X \smallcoprod Y\) by \(\emb_1\). Given morphisms \(f\colon X \map Z\) and \(g\colon Y \map Z\) in \(\mathcal{C}\), the induced morphism \(X \smallcoprod Y \map Z\) is denoted by \(\begin{smallpmatrix} f \\ g \end{smallpmatrix}\).
\item Given a category \(\mathcal{C}\) and an object \(X \in \Ob \mathcal{C}\), the category of objects in \(\mathcal{C}\) under \(X\) will be denoted by \((X \downarrow \mathcal{C})\). The objects in \((X \downarrow \mathcal{C})\) are denoted by \((Y, f)\), where \(Y \in \Ob \mathcal{C}\) and \(f\colon X \map Y\) is a morphism in \(\mathcal{C}\).
\end{itemize}

\section{Preliminaries from homotopical algebra} \label{sec:preliminaries_from_homotopical_algebra}

We recall some basic facts from homotopical algebra. Our main reference is \cite[ch.\ I, \S 1]{quillen:1967:homotopical_algebra}.

\subsection*{Model categories} \label{ssec:model_categories}

Throughout this note, we let \(\mathcal{M}\) be a model category, cf.\ \cite[ch.\ I, \S 1, def.\ 1]{quillen:1967:homotopical_algebra}. In \(\mathcal{M}\), there are three kinds of distinguished morphisms, called \newnotion{cofibrations}, \newnotion{fibrations} and \newnotion{weak equivalences}. Cofibrations are closed under pushouts. If weak equivalences in \(\mathcal{M}\) are closed under pushouts along cofibrations, \(\mathcal{M}\) is said to be \newnotion{left proper}, cf.\ \cite[def.\ 13.1.1(1)]{hirschhorn:2003:model_categories_and_their_localizations}.

An object \(X \in \Ob \mathcal{M}\) is said to be \newnotion{fibrant} if the unique morphism \(\mathcal{M} \map \terminalobject\) is a fibration, where \(\terminalobject\) is a (chosen) terminal object in \(\mathcal{M}\). The full subcategory of \(\mathcal{M}\) of fibrant objects is denoted by \(\Fib(\mathcal{M})\).

The \newnotion{homotopy category} of \(\mathcal{U} \in \{\mathcal{M}, \Fib(\mathcal{M})\}\) is a localisation of \(\mathcal{U}\) with respect to the weak equivalences in \(\mathcal{U}\) and is denoted by \(\HomotopyCategory \mathcal{U}\). The localisation functor of \(\HomotopyCategory \mathcal{U}\) is denoted by \(\LocalisationFunctor = \LocalisationFunctor[\HomotopyCategory \mathcal{U}]\colon \mathcal{U} \map \HomotopyCategory \mathcal{U}\).

Given an object \(X \in \Ob \mathcal{M}\), the category \((X \downarrow \mathcal{M})\) of objects under \(X\) obtains a model category structure where a morphism in \((X \downarrow \mathcal{M})\) is a weak equivalence resp.\ a cofibration resp.\ a fibration if and only if it is one in \(\mathcal{M}\).

\subsection*{Homotopies} \label{ssec:homotopies}

A \newnotion{cylinder} for an object \(X \in \Ob \mathcal{M}\) consists of an object \(Z \in \Ob \mathcal{M}\), a cofibration \(\begin{smallpmatrix} \ins_0 \\ \ins_1 \end{smallpmatrix} = \ins = \ins^Z\colon X \smallcoprod X \map Z\) and a weak equivalence \(\cylinderequivalence = \cylinderequivalence^Z\colon Z \map X\) such that \(\ins \, \cylinderequivalence = \begin{smallpmatrix} 1 \\ 1 \end{smallpmatrix}\). 

Given parallel morphisms \(f, g\colon X \map Y\) in \(\mathcal{M}\), we say that \(f\) is \newnotion{cylinder homotopic} to \(g\), written \(f \cylinderhomotopic g\), if there exists a cylinder \(Z\) for \(X\) and a morphism \(H\colon Z \map Y\) with \(\ins_0 H = f\) and \(\ins_1 H = g\). In this case, \(H\) is said to be a \newnotion{cylinder homotopy} from \(f\) to \(g\). (In the literature, cylinder homotopy is also called left homotopy, cf.\ \cite[ch.\ I, \S 1, def.\ 3, def.\ 4, lem.\ 1]{quillen:1967:homotopical_algebra}.) The relation \(\cylinderhomotopic\) is reflexive and symmetric, but in general not transitive. Moreover, \(\cylinderhomotopic\) is compatible with composition in \(\Fib(\mathcal{M})\). We denote by \(\Fib(\mathcal{M}) / {\cylinderhomotopic}\) the quotient category of \(\Fib(\mathcal{M})\) with respect to the congruence generated by \(\cylinderhomotopic\).

\subsection*{Quillen's homotopy category theorem} \label{ssec:quillens_homotopy_category_theorem}

There are dual notions to fibrant objects, cylinders, cylinder homotopic \(\cylinderhomotopic\), the full subcategory of fibrant objects \(\Fib(\mathcal{M})\), its quotient category \(\Fib(\mathcal{M}) / {\cylinderhomotopic}\) and its homotopy category \(\HomotopyCategory \Fib(\mathcal{M})\), namely \newnotion{cofibrant objects}, \newnotion{path objects}, \newnotion{path homotopic} \(\pathhomotopic\), the full subcategory of cofibrant objects \(\Cof(\mathcal{M})\), its quotient category \(\Cof(\mathcal{M}) / {\pathhomotopic}\) and its homotopy category \(\HomotopyCategory \Cof(\mathcal{M})\), respectively. Moreover, an object \(X \in \Ob \mathcal{M}\) is said to be bifibrant if it is cofibrant and fibrant. On the full subcategory of bifibrant objects \(\Bif(\mathcal{M})\), the relations \(\cylinderhomotopic\) and \(\pathhomotopic\) coincide and yield a congruence. One writes \({\homotopic} := {\cylinderhomotopic} = {\pathhomotopic}\) in this case, and the quotient category is denoted by \(\Bif(\mathcal{M}) / {\homotopic}\). Moreover, \(\HomotopyCategory \Bif(\mathcal{M})\) is a localisation of \(\Bif(\mathcal{M})\) with respect to the weak equivalences in \(\Bif(\mathcal{M})\).

Quillen's homotopy category theorem \cite[ch.\ I, \S 1, th.\ 1]{quillen:1967:homotopical_algebra} (cf.\ also \cite[cor.\ 1.2.9, th.\ 1.2.10]{hovey:1999:model_categories}) states that the various inclusion and localisation functors induce the following commutative diagram, where the functors labeled by \(\categoricalequivalent\) are equivalences and the functor labeled by \(\isomorphic\) is an isofunctor.
\[\begin{tikzpicture}[baseline=(m-3-1.base)]
  \matrix (m) [matrix of math nodes, row sep=2.5em, column sep=2.5em, text height=1.6ex, text depth=0.45ex, inner sep=0pt, nodes={inner sep=0.3333em}]{
    \Cof(\mathcal{M}) / {\pathhomotopic} & \HomotopyCategory \Cof(\mathcal{M}) & \\
    \Bif(\mathcal{M}) / {\homotopic} & \HomotopyCategory \Bif(\mathcal{M}) & \HomotopyCategory \mathcal{M} \\
    \Fib(\mathcal{M}) / {\cylinderhomotopic} & \HomotopyCategory \Fib(\mathcal{M}) & \\};
  \path[->, font=\scriptsize]
    (m-1-1) edge (m-1-2)
    (m-1-2) edge node[right=1pt] {\(\categoricalequivalent\)} (m-2-3)
    (m-2-1) edge node[above] {\(\isomorphic\)} (m-2-2)
            edge (m-3-1)
            edge (m-1-1)
    (m-2-2) edge node[right] {\(\categoricalequivalent\)} (m-3-2)
            edge node[right] {\(\categoricalequivalent\)} (m-1-2)
    (m-3-1) edge (m-3-2)
    (m-3-2) edge node[right=1pt] {\(\categoricalequivalent\)} (m-2-3);
\end{tikzpicture}\]
In this note, we treat the question whether the functors \(\Fib(\mathcal{M}) / {\cylinderhomotopic} \map \HomotopyCategory \Fib(\mathcal{M})\) and \(\Cof(\mathcal{M}) / {\pathhomotopic} \map \HomotopyCategory \Cof(\mathcal{M})\) are faithful. By duality, it suffices to consider the first functor.

\subsection*{The model category \(\fgMod(\Integers / 4)\)} \label{ssec:the_model_category_modz/4}

The category \(\fgMod(\Integers / 4)\) of finitely generated modules over \(\Integers / 4\) is a Frobenius category (with respect to all short exact sequences), that is, there are enough projective and injective objects in \(\fgMod(\Integers / 4)\) and, moreover, these objects coincide (we call such objects bijective). Therefore \(\fgMod(\Integers / 4)\) carries a canonical model category structure (cf.\ also \cite[sec.\ 2.2]{hovey:1999:model_categories}): The cofibrations are the monomorphisms and the fibrations are the epimorphisms in \(\fgMod(\Integers / 4)\). Every object in \(\fgMod(\Integers / 4)\) is bifibrant, and the weak equivalences are precisely the homotopy equivalences, where parallel morphisms \(f\) and \(g\) are homotopic if \(g - f\) factors over a bijective object in \(\fgMod(\Integers / 4)\). That is, the weak equivalences in \(\fgMod(\Integers / 4)\) are the stable isomorphisms and the homotopy category of \(\fgMod(\Integers / 4)\) is isomorphic to the stable category of \(\fgMod(\Integers / 4)\), cf.\ \cite[ch.\ I, sec.\ 2.2]{happel:1988:triangulated_categories_in_the_representation_theory_of_finite-dimensional_algebras}.

We remark that every object in \(\fgMod(\Integers / 4)\) is isomorphic to \((\Integers / 4)^{\directsum k} \directsum (\Integers / 2)^{\directsum l}\) for some \(k, l \in \Naturals_0\), and every bijective object is isomorphic to \((\Integers / 4)^{\directsum k}\) for some \(k \in \Naturals_0\).

\section{\texorpdfstring{Faithfulness of the functor \(\Fib(\mathcal{M}) / {\cylinderhomotopic} \map \HomotopyCategory \Fib(\mathcal{M})\)}{Faithfulness of the functor Fib(M) / cylinderhomotopic -> Ho Fib(M)}} \label{sec:faithfulness_of_the_functor_fib/ch_to_ho_fib}

We give a sufficient criterion for the functor under consideration to be faithful.

\begin{proposition*}
If the model category \(\mathcal{M}\) is left proper and if \(w \smallcoprod w\) is a weak equivalence for every weak equivalence \(w\) in \(\mathcal{M}\), then \(\cylinderhomotopic\) is a congruence on \(\Fib(\mathcal{M})\) and the canonical functor \(\Fib(\mathcal{M}) / {\cylinderhomotopic} \map \HomotopyCategory \Fib(\mathcal{M})\) is faithful.
\end{proposition*}
\begin{proof}
We suppose given fibrant objects \(X\) and \(Y\) and morphisms \(f, g\colon X \map Y\) with \(\LocalisationFunctor f = \LocalisationFunctor g\) in \(\HomotopyCategory \Fib(\mathcal{M})\). By \cite[th.\ 1(ii)]{brown:1974:abstract_homotopy_theory_and_generalized_sheaf_cohomology}, there exists a weak equivalence \(w\colon X' \map X\) such that \(w f \pathhomotopic w g\). It follows that \(w f \cylinderhomotopic w g\) by \cite[ch.\ I, \S 1, dual of lem.\ 5]{quillen:1967:homotopical_algebra}, that is, there exists a cylinder \(Z'\) for \(X'\) and a cylinder homotopy \(H'\colon Z' \map Y\) from \(w f\) to \(w g\). We let
\[\begin{tikzpicture}[baseline=(m-2-1.base)]
  \matrix (m) [diagram]{
    X' \smallcoprod X' & X \smallcoprod X \\
    Z' & Z \\};
  \path[->, font=\scriptsize]
    (m-1-1) edge node[above] {\(w \smallcoprod w\)} node[sloped, below] {\(\weaklyequivalent\)} (m-1-2)
            edge node[left] {\(\ins^{Z'}\)} (m-2-1)
    (m-1-2) edge node[right] {\(i\)} (m-2-2)
    (m-2-1) edge node[above] {\(w'\)} node[sloped, below] {\(\weaklyequivalent\)} (m-2-2);
\end{tikzpicture}\]
be a pushout of \(w \smallcoprod w\) along \(\ins^{Z'}\). By assumption, \(w \smallcoprod w\) and \(w'\) are weak equivalences. Since \((w \smallcoprod w) \begin{smallpmatrix} 1 \\ 1 \end{smallpmatrix} = \ins^{Z'} \cylinderequivalence^{Z'} w\), there exists a unique morphism \(s\colon Z \map X\) with \(\begin{smallpmatrix} 1 \\ 1 \end{smallpmatrix} = i s\) and \(\cylinderequivalence^{Z'} w = w' s\). Then \(s\) is a weak equivalence since \(\cylinderequivalence^{Z'}\), \(w\) and \(w'\) are weak equivalences and therefore \(Z\) becomes a cylinder for \(X\) with \(\ins^{Z} := i\) and \(\cylinderequivalence^{Z} := s\).
Moreover, \((w \smallcoprod w) \begin{smallpmatrix} f \\ g \end{smallpmatrix} = \ins^{Z'} H'\) implies that there exists a unique morphism \(H\colon Z \map Y\) with \(\begin{smallpmatrix} f \\ g \end{smallpmatrix} = \ins^Z H\) and \(H' = w' H\). So in particular \(f \cylinderhomotopic g\). 
\[\begin{tikzpicture}[baseline=(m-3-1.base)]
  \matrix (m) [diagram]{
    X' \smallcoprod X' & X \smallcoprod X & Y \\
    Z' & Z & Y \\
    X' & X & \\};
  \path[->, font=\scriptsize]
    (m-1-1) edge node[above] {\(w \smallcoprod w\)} node[sloped, below] {\(\weaklyequivalent\)} (m-1-2)
            edge node[left] {\(\ins^{Z'}\)} (m-2-1)
    (m-1-2) edge node[above] {\(\begin{smallpmatrix} f \\ g \end{smallpmatrix}\)} (m-1-3)
            edge node[right] {\(\ins^Z\)} (m-2-2)
    (m-1-3) edge[equality] (m-2-3)
    (m-2-1) edge node[above] {\(w'\)} node[sloped, below] {\(\weaklyequivalent\)} (m-2-2)
            edge node[left] {\(\cylinderequivalence^{Z'}\)} node[sloped, above] {\(\weaklyequivalent\)} (m-3-1)
    (m-2-2) edge node[above] {\(H\)} (m-2-3)
            edge node[right] {\(\cylinderequivalence^{Z}\)} node[sloped, below] {\(\weaklyequivalent\)} (m-3-2)
    (m-3-1) edge node[above] {\(w\)} node[sloped, below] {\(\weaklyequivalent\)} (m-3-2);
\end{tikzpicture}\]

Altogether, we have shown that morphisms in \(\Fib(\mathcal{M})\) represent the same morphism in \(\HomotopyCategory \Fib(\mathcal{M})\) if and only if they are cylinder homotopic. In particular, \(\cylinderhomotopic\) is a congruence on \(\Fib(\mathcal{M})\).
\end{proof}

The following counterexample shows that the canonical functor \(\Fib(\mathcal{M}) / {\cylinderhomotopic} \map \HomotopyCategory \Fib(\mathcal{M})\) is not faithful in general.

\begin{example*}
We consider the category \((\Integers / 4 \downarrow \fgMod(\Integers / 4))\) of finitely generated \(\Integers / 4\)-modules under \(\Integers / 4\) with the model category structure inherited from \(\fgMod(\Integers / 4)\), cf.\ \ref{sec:preliminaries_from_homotopical_algebra}. All objects of \((\Integers / 4 \downarrow \fgMod(\Integers / 4))\) are fibrant since all objects in \(\fgMod(\Integers / 4)\) are fibrant.

We study morphisms \((\Integers / 4, 2) \map (\Integers / 4 \directsum \Integers / 2, \begin{smallpmatrix} 2 & 0 \end{smallpmatrix})\) in \((\Integers / 4 \downarrow \fgMod(\Integers / 4))\). We let \((Z, t)\) be a cylinder of \((\Integers / 4, 2)\) and we let \(H\colon (Z, t) \map (\Integers / 4 \directsum \Integers / 2, \begin{smallpmatrix} 2 & 0 \end{smallpmatrix})\) be a cylinder homotopy (from \(\ins_0 H\) to \(\ins_1 H\)). Then we have a weak equivalence \((Z, t) \map (\Integers / 4, 2)\) in \((\Integers / 4 \downarrow \fgMod(\Integers / 4))\) and hence a weak equivalence \(Z \map \Integers / 4\) in \(\fgMod(\Integers / 4)\). Thus \(Z\) is bijective and therefore we may assume that \(Z = (\Integers / 4)^{\directsum k}\). Since \(\ins_0\) and \(\ins_1\) are morphisms from \((\Integers / 4, 2)\) to \((Z, t)\), we have \(2 \ins_0 = t = 2 \ins_1\) and hence \(\ins_0 \congruent_2 \ins_1\) as morphisms from \(\Integers / 4\) to \(Z\). But this implies that the second components of \(\ins_0 H\) and \(\ins_1 H\) are the same. In other words, we have shown that cylinder homotopic morphisms from \((\Integers / 4, 2)\) to \((\Integers / 4 \directsum \Integers / 2, \begin{smallpmatrix} 2 & 0 \end{smallpmatrix})\) coincide in the second component. It follows that the morphisms \(\begin{smallpmatrix} 1 & 0 \end{smallpmatrix}\colon (\Integers / 4, 2) \map (\Integers / 4 \directsum \Integers / 2, \begin{smallpmatrix} 2 & 0 \end{smallpmatrix})\) and \(\begin{smallpmatrix} 1 & 1 \end{smallpmatrix}\colon (\Integers / 4, 2) \map (\Integers / 4 \directsum \Integers / 2, \begin{smallpmatrix} 2 & 0 \end{smallpmatrix})\) in \((\Integers / 4 \downarrow \fgMod(\Integers / 4))\) represent different morphisms in the quotient category \(\Fib((\Integers / 4 \downarrow \fgMod(\Integers / 4))) / {\cylinderhomotopic}\). 

On the other hand, since \(\Integers / 4\) is bijective, the morphism \(2\colon \Integers / 4 \map \Integers / 4\) is a weak equivalence in \(\fgMod(\Integers / 4)\), and therefore \(2\colon (\Integers / 4, 1) \map (\Integers / 4, 2)\) is a weak equivalence in \((\Integers / 4 \downarrow \fgMod(\Integers / 4))\). But \(2 \begin{smallpmatrix} 1 & 0 \end{smallpmatrix} = 2 \begin{smallpmatrix} 1 & 1 \end{smallpmatrix}\) as morphisms from 
\((\Integers / 4, 1)\) to \((\Integers / 4 \directsum \Integers / 2, \begin{smallpmatrix} 2 & 0 \end{smallpmatrix})\) in \((\Integers / 4 \downarrow \fgMod(\Integers / 4))\), so in particular \(\LocalisationFunctor(2 \begin{smallpmatrix} 1 & 0 \end{smallpmatrix}) = \LocalisationFunctor(2 \begin{smallpmatrix} 1 & 1 \end{smallpmatrix})\) and hence \(\LocalisationFunctor \begin{smallpmatrix} 1 & 0 \end{smallpmatrix} = \LocalisationFunctor \begin{smallpmatrix} 1 & 1 \end{smallpmatrix}\).
\end{example*}

% bibliography

\bigskip

{\raggedleft William G.\ Dwyer \\ Department of Mathematics \\ University of Notre Dame \\ Notre Dame, IN 46556 \\ dwyer.1@nd.edu \\ \url{http://www.nd.edu/~wgd/} \\}

\bigskip

{\raggedleft Andrei R\u{a}dulescu-Banu \\ 86 Cedar St \\ Lexington, MA 02421 USA \\ andrei@alum.mit.edu \\}

\bigskip

{\raggedleft Sebastian Thomas \\ Lehrstuhl D für Mathematik \\ RWTH Aachen \\ Templergraben 64 \\ D-52062 Aachen \\ sebastian.thomas@math.rwth-aachen.de \\ \url{http://www.math.rwth-aachen.de/~Sebastian.Thomas/} \\}

\end{document}